\documentclass[12pt, a4paper]{amsart}

\usepackage{amssymb}
\usepackage{stmaryrd}
\usepackage[dvips]{graphics}

\theoremstyle{plain}
\newtheorem*{thm}{Theorem}
\newtheorem*{lem}{Lemma}

\newtheorem{theorem}{Theorem}[section]
\newtheorem{lemma}[theorem]{Lemma}
\newtheorem{corollary}[theorem]{Corollary}
\newtheorem{proposition}[theorem]{Proposition}

\theoremstyle{definition}
\newtheorem{definition}[theorem]{Definition}

\newtheorem*{qu}{Question}

\newtheorem{remark}[theorem]{Remark}

\theoremstyle{remark}
\newtheorem{clm}{Claim}
\newtheorem*{claim}{Claim}
\newtheorem*{acknowledgements}{Acknowledgements}

\numberwithin{figure}{section}           
\numberwithin{table}{section}
\numberwithin{equation}{section}

\DeclareMathOperator{\aut}{Aut}
\DeclareMathOperator{\autf}{Aut_{\it f}}

\DeclareMathOperator{\Som}{V}
\DeclareMathOperator{\Are}{E}
\DeclareMathOperator{\Ch}{Chain}
\DeclareMathOperator{\B}{End}

\DeclareMathOperator{\dist}{dist}

\DeclareMathOperator{\stab}{Stab}
\DeclareMathOperator{\fix}{Fix}
\DeclareMathOperator{\im}{Im}

\DeclareMathOperator{\PSL}{PSL}

\newcommand{\N}{{\mathbb N}}
\newcommand{\Z}{{\mathbb Z}}

\newcommand{\R}{{\mathbb R}}

\newcommand{\U}{{\mathcal U}}

\newcommand{\K}{{\mathcal K}}

\newcommand{\G}{{\mathcal G}}

\newcommand{\card}{{\rm Card}}
\newcommand{\id}{{\rm id}}
\newcommand{\ram}{{\rm ram}}

\begin{document}

\title{Boundedly simple groups of automorphisms of trees}
\author{Jakub Gismatullin}
\thanks{The author is supported by Polish Goverment MNiSW grants N N201 384134, N N201 545938, by fellowship START of the Foundation for Polish Science and by the Marie Curie Intra-European Fellowship MODGROUP no. PIEF-GA-2009-254123.}
\date{\today}
\address{Instytut Matematyczny Uniwersytetu Wroc{\l}awskiego, pl. Grunwaldzki 2/4, 50-384 Wroc{\l}aw, Poland}
\address{and}
\address{Institute of Mathematics of the Polish Academy of Sciences, ul. \'Sniadeckich 8, 00-956 Warsaw, Poland}
\address{and}
\address{School of Mathematics, University of Leeds, Woodhouse Lane, Leeds, LS2 9JT, UK}
\email{gismat@math.uni.wroc.pl, j.gismatullin@leeds.ac.uk}


\begin{abstract}
A group is boundedly simple if, for some constant $N$, every nontrivial conjugacy class generates the whole group in $N$ steps. For a large class of trees, Tits proved simplicity of a canonical subgroup of the automorphism group, which is generated by pointwise stabilizers of edges. We prove that only for uniform subdivisions of biregular trees are such groups boundedly simple. In fact these groups are 8-boundedly simple. As a consequence, we prove that if $G$ is boundedly simple (or from a certain class $\K$) and $G$ acts by automorphisms on a tree, then $G$ fixes some vertex of $A$, or stabilizes some end of $A$, or the smallest nonempty $G$-invariant subtree of $A$ is a uniform subdivision of a biregular tree.
\end{abstract}

\maketitle

\section{Introduction}

A group $G$ is simple (in the algebraic sense) if and only if $G$ is generated by every nontrivial conjugacy class. A finer notion is that of \emph{bounded simplicity}. A group $G$ is called \emph{$N$-boundedly simple} if for every two nontrivial elements $g,h\in G$, the element $h$ is the product of $N$ or fewer conjugates of $g^{\pm1}$, i.e. \[G = \left(g^G \cup {g^{-1}}^G\right)^{\leq N}.\] We say $G$ is \emph{boundedly simple} if it is $N$-boundedly simple, for some $N\in\N$. 

In this paper we are interested in actions of boundedly simple groups on trees. Our results were inspired by the following theorem due to Tits.

\begin{thm}{\cite[Theorem 4.5]{tits}}
Suppose that $A$ is a tree and $G$ is a group acting by automorphisms on $A$ without leaving invariant any nonempty proper subtree of $A$ or any end of $A$. Assume that $G$ has Tits' independence property $(P)$ (see Definition \ref{def:p}). Let $G^+$ be the subgroup of $G$ generated by pointwise stabilizers in $G$ of edges of $A$ (see Definition \ref{def:ram}$(2)$). Then $G^+$ is a simple group. Furthermore, every subgroup of $G$ normalized by $G^+$ is trivial or contains $G^+$.
\end{thm}

The full group of automorphisms $\aut(A)$ has property $(P)$ and in many cases does not leave invariant subtrees or ends of $A$. In such a case, by the above theorem, $\aut^+(A)$ is simple. We determine trees such that $\aut^+(A)$ is boundedly simple. In fact, we consider a more general situation of a tree with a coloring $f$ of the set $\Som(A)$ of all vertices and group $\autf^+(A)$ of color-preserving automorphisms.

By $A_{n,m}$, for some cardinal numbers $m,n\geq 3$, we denote an \emph{$(n,m)$-regular (biregular) tree}. That is, a tree in which every vertex is black or white with vertices of the same color non-adjacent, every white vertex is connected with $n$ black vertices and every black vertex is connected with $m$ white vertices. By an \emph{$m$-regular tree} we mean $A_{m,m}$. We prove the following structure theorem of the automorphism group of a biregular tree.

\begin{thm}{\bf \ref{thm:bireg}}
Suppose that $n, m\geq 3$ are cardinals and $A_{n,m}$ is an $(n,m)$-regular (biregular) tree. The group $\aut^+(A_{n,m})$ is $8$-boundedly simple. Moreover, if $m=n$, then $\left[\aut(A_{n,m}):\aut^+(A_{n,m})\right]= 2$; if $m\ne n$, then $\aut(A_{n,m})^+=\aut(A_{n,m})$.
\end{thm}

By a \emph{uniform subdivision of a tree} we mean roughly a subdividing of each edge of the tree into the same number of edges (see Definition \ref{def:almost}).

\begin{thm}{\bf \ref{thm:main}}
Assume that $(A, f\colon \Som(A) \to I)$ is a colored tree and $\autf^+(A)$ is boundedly simple and nontrivial. Then $\autf^+(A)$ fixes some vertex of $A$, or leaves invariant some end of $A$, or leaves invariant a subtree $A' \subseteq A$, which is a uniform subdivision of $(n,m)$-regular tree, for some $n,m\geq 3$.

In particular, if $\autf^+(A)$ leaves no nonempty proper subtree of $A$ invariant and stabilizes no end, then $A$ is a uniform subdivision of a biregular tree and $\autf^+(A)$ is $8$-boundedly simple.
\end{thm}

As a consequence, the bounded simplicity of automorphism groups characterises the biregular trees. We do not expect that the bound 8 is sharp. The proof of Theorem \ref{thm:main} goes through Proposition \ref{prop:autf}, which asserts that if $\autf^+(A)$ is boundedly simple, then some configuration in the code of $A$ is forbidden.

In the last section we deal with a more general set-up of an action of a group on a tree. We consider groups from a certain class $\K$ (see Definition \ref{def:k}), consisting of all groups $G$ such that $G$ and all subgroups of index 2 of $G$ are boundedly generated by some finite set of conjugacy classes. In particular $\K$ contains all boundedly simple groups. Our motivation for studying such actions comes from Bruhat-Tits buildings for $\PSL_2(K)$, where $K$ is a field with a discrete valuation (see \cite[Chapter II]{serre}). That is, $\PSL_2(K)$ acts by automorphisms on an $(n+1)$-regular tree $A_{n+1,n+1}$ (its Bruhat-Tits building), where $n$ is the cardinality of the residue field. In fact, $\PSL_2(K)$ is a subgroup of $\aut^+(A_{n+1,n+1})$, and leaves no nonempty proper subtree of $A_{n+1,n+1}$ invariant and does not stabilize any end of $A_{n+1,n+1}$. On the other hand, it is well known that for an arbitrary field $K$, the group $\PSL_2(K)$ is boundedly simple (by \cite{psl}, $\PSL_2(K)$ is $5$-boundedly simple), so $\PSL_2(K)$ is in class $\K$.

\begin{thm}{\bf \ref{thm:act}}
Suppose that $A$ is a tree, $G<\aut(A)$ and $G$ is from the class $\K$. If $G$ leaves no nonempty proper subtree of $A$ invariant and does not stabilize any end of $A$, then $A$ is a uniform subdivision of some $(n,m)$-regular tree, for some $n,m\geq 3$.
\end{thm}

As an immediate consequence of the above theorem, we have the following `Invariant subtree theorem for boundedly simple groups' (Corollary \ref{cor:act}): 
\begin{quote}
If $G<\aut(A)$ and $G\in\K$, then $G$ fixes some vertex of $A$, or stabilizes some end of $A$, or the smallest nonempty $G$-invariant subtree of $A$ is a uniform subdivision of a biregular tree.
\end{quote}

Bounded simplicity arises naturally in model theory in the study of first order expressibility of simplicity for groups. For fixed $N$, the property of `$N$-bounded simplicity' is first order expressible, i.e. can be written as a sentence in the first order logic. Therefore, for each $N\in\N$, the class of $N$-boundedly simple groups is an \emph{elementary class} (or an \emph{axiomatizable class}) of structures. Every elementary class of structures is closed under taking ultraproducts (and elementary extensions). Some boundedly simple groups has been constructed in \cite{muranov}. In fact, the following well known lemma characterises bounded simplicity.

\begin{lem} \label{lem}
The following conditions are equivalent for any group $G$.
\begin{enumerate}
\item[(1)] $G$ is boundedly simple.
\item[(2)] Some ultrapower $G^{\N}/\U$ of $G$ over some non-principal ultrafilter $\U$ is a simple group.
\end{enumerate}
\end{lem}
\begin{proof}
$(1)\Rightarrow(2)$ Bounded simplicity is a first order property and, by {\L}o{\'s} Theorem, ultrapowers preserve first-order conditions. Hence every ultrapower is boundedly simple, and thus simple.

$(2)\Rightarrow(1)$ Let an ultrapower $G^{\N}/\U$ be simple. Assume contrary to $(1)$, that for every $N\in\N$, there is $g_N \in G\setminus \{e\}$ and \[h_N \in G \setminus \left(g_N^G \cup {g_N^{-1}}^G\right)^{\leq N}.\] Consider $g=(g_N)_{N\in\N}/\U$ and $h=(h_N)_{N\in\N}/\U$ from $G^{\N}/\U$. Then the normal closure \[ H = \bigcup_{n<\N} \left(g^{G^{\N}/\U} \cup {g^{-1}}^{G^{\N}/\U}\right)^n\] of $g$ in $G^{\N}/\U$ is a nontrivial subgroup of $G^{\N}/\U$, which is proper (as $h\not\in H$); this is impossible.
\end{proof}

Using the above lemma one can give an easy proof of bounded simplicity of $\PSL_n(K)$, where $n\geq2$ and $|K|\geq 4$. Namely, let $K$ be an arbitrary field with $|K|\geq 4$. Then $\PSL_n(K)^{\N}/\U \cong \PSL_n\left(K^{\N}/\U\right)$. However $\PSL_n(F)$ is a simple group, for an arbitrary field $F$ with $|F|\geq 4$. Hence by the lemma $\PSL_n(K)$ is boundedly simple. In fact, by \cite[Theorem M]{gord}, $G(k)$ is boundedly simple, where $k$ is a field, and $G$ is any $k$-split, semisimple, simply connected linear algebraic group (that is a Chevalley group).

There are many fixed point results for actions of linear groups on tress, or on some other spaces. For example, the following fact is due to Tits \cite[Corollary 4]{tits_lie}. Suppose $G$ is an almost simple isotropic linear algebraic $k$-group and let $G(k)^+$ be a Zariski dense subgroup of $G(k)$ generated by rational unipotent elements. Assume that the $k$-rank of $G$ is at least 2 and $G(k)^+$ acts on an $\R$-tree $A$ by isometries. Then $G(k)^+$ fixes some vertex of $A$ or centralizes an end of $A$.

Serre introduced groups with the $(FA)$ property. A group $G$ has $(FA)$ if any action of $G$ on a tree without inversion has a fixed point. Recently in \cite{farb} fixed point theorems of action of certain groups on nonpositively curved space have been proven, and stronger notions than $(FA)$ were considered.  However, certain group from the class $\K$ do not have $(FA)$, for example $\PSL_2(K)$ for some field $K$, so our results (Theorem \ref{thm:act}, Corollary \ref{cor:act}) for these groups are new.





Suppose $G$ is any group acting on an $\R$-tree. In \cite{tignol}, the simplicity of the subgroup $G^+$ of $G$ has been obtained under similar assumptions as in \cite[Theorem 4.5]{tits}. In \cite{hag}, the authors study groups of automorphisms of some negatively curved spaces (hyperbolic buildings, Cayley graphs of word hyperbolic Coxeter groups and generalised cubical complexes). In particular they show that the group of type-preserving automorphisms of $I_{p,q}$ (for $p\geq 5$, $q\geq 3$) is simple; here $I_{p,q}$ denotes the $2$-dimensional hyperbolic building of M. Bourdon \cite{bur}. It would be interesting to generalise results from this paper to spaces considered in \cite{tignol, hag, bur}.

\begin{qu}
Are the simple groups considered in \cite{tignol, hag} boundedly simple? 
\end{qu}

\section{Basic Notation and Prerequisites}

We use the notation and basic facts from \cite{tits}. A \emph{tree} is a connected graph without cycles. In this paper $A$ always denotes a tree. By $\Som(A)$ we denote the set of vertices of $A$. The set of edges $\Are(A)$ is a collection of some 2-element subsets of $\Som(A)$. A sequence of vertices $(s_i)_{-m < i < n, i\in\Z}$, where $n,m\in\N\cup\{\omega\}$, is a \emph{chain in $A$} if
\begin{itemize}
\item the vertices $s_i$, for $-m < i < n$, are all distinct,
\item if $n+m>2$, then for each $-m < i < n-1$, $\{s_i,s_{i+1}\}$ is an edge of $A$.
\end{itemize}
Suppose $C=(s_i)_{-m < i < n}$ is a chain. If $n,m\in\N$ and $n+m>2$, then we also call $C$ a \emph{chain of length $n+m-2$ joining $s_{-m+1}$ and $s_{n-1}$}. If $n=\omega$ and $m\in\N$, then we call $C$ a \emph{one-way infinite chain from $s_{-m+1}$} (similarly, when $m=\omega$ and $n\in\N$). When $n=m=\omega$, we call $C$ a \emph{two-way infinite chain}. Note that for any two distinct
vertices of a tree, there exists a unique chain joining them. Let $\Ch(A)$ be the set of all one-way infinite chains starting at some vertex of $A$. \emph{Ends} are equivalence classes of the following relation defined on $\Ch(A)$: $C \sim C' \Leftrightarrow C \cap C' \in \Ch(A)$. The set of ends is denoted by $\B(A)$.

By $\aut(A)$ we denote the group of all automorphisms of $A$, i.e. permutations of $\Som(A)$ preserving edges. An automorphism $\alpha \in \aut(A)$ is called a \emph{rotation} if it stabilizes some vertex $s\in \Som(A)$, i.e. $\alpha(s) = s$. We say $\alpha$ is an \emph{inversion} if for some edge $\left\{s,s'\right\}\in\Are(A)$, $\alpha(s)=s'$ and $\alpha\left(s'\right) = s$. If for some two-way infinite chain $C$ in $A$, an automorphism $\alpha$ leaves $C$ invariant and is not a rotation or an inversion, then we call $\alpha$ a \emph{translation}; in this case $C$ is the unique two-way infinite chain with the above properties and $\alpha$ restricted to $C$ is a nontrivial translation. We also call $C$ the \emph{axis of $\alpha$}. The \emph{translation length} of $\alpha$ is the infimum of the distances between $s$ and $\alpha(s)$, for all $s\in\Som(A)$. Note that the translation length of an arbitrary translation is always positive. By \cite[Proposition 3.2]{tits} the group $\aut(A)$ is a disjoint union of rotations, inversions and translations. The subtree of $A$ consisting of vertices fixed pointwise by $\alpha$ is called the \emph{fixed tree} of $\alpha$ and is denoted by $\fix(\alpha)$. The subgroup of $\aut(A)$ stabilizing pointwise a given subtree $A'$ of $A$ is denoted by $\stab\left(A'\right)$. For $G<\aut(A)$, by $\stab^G\left(A'\right)$ we denote $\stab\left(A'\right)\cap G$. The group $\aut(A)$ acts naturally on the set $\B(A)$ of ends of $A$. 

\begin{definition}{\cite[2.5]{tits}} \label{def:centr} Let $\alpha\in \aut(A)$ and $b\in \B(A)$. We say that
\begin{itemize}
\item[(1)] $\alpha$ \emph{stabilizes} $b$ or $\alpha$ \emph{leaves invariant} $b$, if $\alpha(b) = b$;
\item[(2)] $\alpha$ \emph{centralizes} $b$, if $\alpha$ fixes pointwise some chain $C$ from $b$. The set of all elements that centralize $b$ forms a group, called the \emph{centralizer} of $b$.
\end{itemize}
\end{definition}

Clearly, if $\alpha$ centralizes $b$, then $\alpha$ also stabilizes $b$. If $\alpha$ is not a nontrivial translation, then the converse is also true, i.e. if $\alpha(b) = b$, then for some $C\in b$, $\alpha_{|C}=\id_C$. To see this, note that $\alpha$ must be a rotation, i.e. $\alpha(s)=s$ for some $s$. Then $\alpha$ fixes pointwise some infinite chain $C$ from $b$ starting at $s$.

The next two lemmas are well known (see e.g. \cite[Section 6.5]{serre}). However, for the completeness of the exposition we provide proofs.

\begin{lemma} \label{lem:trans}
%
Suppose that $\alpha, \beta\in\aut(A)$ are rotations and $\fix(\alpha)\cap\fix(\beta)=\emptyset$. Then $\alpha\circ \beta$ is a translation with an even translation length.
\end{lemma}
\begin{figure}
\centering
\includegraphics{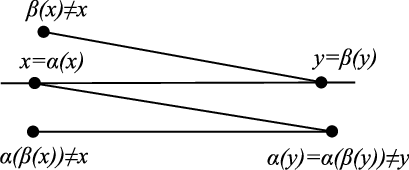}
\caption{Composition of two rotations}
\label{f:trans}
\end{figure}
\begin{proof} Let $\gamma = \alpha\circ\beta$. We use the following criterion \cite[Lemma 3.1]{tits} for an automorphism $\gamma\in\aut(A)$ to be a translation:
\begin{itemize}
\item[$(\spadesuit)$] if for some vertices $x\neq y\in \Som(A)$, $x$ is on the chain joining $y$ with $\gamma(y)$ and $\gamma(y)$ is on the chain joining $x$ with $\gamma(x)$, then $\gamma$ is a translation along a two-way infinite chain containing $y, x, \gamma(y)$ and $\gamma(x)$, with the translation length $\dist(x,\gamma(x)) = \dist(y,\gamma(y))$.
\end{itemize}
One can find $x\neq y\in \Som(A)$ such that $\alpha(x)=x$, $\beta(y)=y$. Suppose $C$ is the chain joining $x$ with $y$. We may assume that the only vertex on $C$ fixed by $\alpha$ is $x$ and the only vertex on $C$ fixed by $\beta$ is $y$.

Since $\alpha(x)=x$, $\alpha(y) = \gamma(y) \neq y$ and $\gamma(x) = \alpha(\beta(x)) \neq \alpha(x) = x$, the chain from $y$ to $\gamma(x)$ first goes through $x$ and then through $\gamma(y)$ (see Figure \ref{f:trans}). Therefore by $(\spadesuit)$, the translation length of $\gamma$ is $\dist(y,\gamma(y)) = 2\dist(y,x)$.
\end{proof}

It is proved in \cite[Proposition 3.4]{tits} that if a subgroup $G < \aut(A)$ does not contain translations, then $G$ pointwise stabilizes some vertex or edge of $A$, or centralizes some end of $A$. The proof of this fact uses the assumption that $G$ is a group in a very limited way, so a slightly more general fact is true (Lemma \ref{lem:stab} below). We use this generalization in the proof of Proposition \ref{prop:tits}.

\begin{lemma} \label{lem:stab}
If $X\subseteq \aut(A)$ and $X \cup XX$ does not contain translations, then the group $\langle X \rangle$ generated by $X$ also does not contain translations. Hence, $X$ fixes some vertex of $A$, or leaves invariant some edge of $A$ or centralizes some end of $A$.
\end{lemma}
\begin{proof} It is enough to prove that $G = \langle X \rangle$ does not contain translations. The rest follows from \cite[Proposition 3.4]{tits}. 

Upon replacing the tree $A$ by its first barycentric subdivision, there is no loss of generality in assuming that $X$ contains no inversions. Hence, every element of $X$ is a rotation. The family $\{\fix(\alpha) : \alpha\in X\}$ has the following property: for every $\alpha, \beta \in X$ \[\fix(\alpha) \cap \fix(\beta) \neq \emptyset; \tag{1}\] for otherwise, by Lemma \ref{lem:trans}, $XX$ contains a translation (note that this is a tree version of Helly's theorem on convex sets). Take arbitrary $\alpha_1,\ldots,\alpha_n$ from $X$ and let $A_i=\fix(\alpha_i)$. We show by induction that $B_1\cap\ldots\cap B_n\neq\emptyset$, whenever $B_1,\ldots,B_n$ are subtrees of $A$ satisfying Helly's condition $(1)$, and hence $A_1\cap\ldots\cap A_n\neq\emptyset$, so $\prod_{1\leq i\leq n}\alpha_i$ is a rotation. For $n=3$, let $s_r \in B_s\cap B_t$, for pairwise distinct $r,s,t$ from $\{1,2,3\}$. If $s$ is the center of the triangle with vertices $s_1$, $s_2$ and $s_3$, then $s\in B_1\cap B_2\cap B_3$. In the general case of $n+1$ subtrees, consider $C_i = B_i\cap B_{n+1}$, where $1\leq i \leq n$. Then $C_i\cap C_j = B_i\cap B_j\cap B_{n+1}\neq\emptyset$, by the case $n=3$. Hence by induction, $\bigcap_{1\leq i \leq n}C_i = \bigcap_{1\leq j \leq n+1}B_j\neq\emptyset$.
\end{proof}

We will deal with some groups of automorphisms of trees which satisfy \emph{Tits' independence property $(P)$} (\cite[4.2]{tits}). Let $G < \aut(A)$ and $C$ be an arbitrary (finite or infinite) chain in $A$. Let $\pi \colon \Som(A) \to \Som(C)$ be the natural projection, so that $\pi(x)\in \Som(C)$ is the closest vertex to $x$. For every $s\in\Som(C)$ there is an induced projection of the stabilizer \[\rho_s \colon \stab^G(C) \longrightarrow \aut(\pi^{-1}[s]).\]

\begin{definition} \label{def:p} We say that $G<\aut(A)$ has the property $(P)$ if for every chain $C$ in $A$, the mapping \[\rho = (\rho_s)_{s\in\Som(C)} \colon \stab^G(C) \longrightarrow \prod_{s\in\Som(C)} \im(\rho_s)\] is an isomorphism.
\end{definition}

For example, the full group of automorphisms $\aut(A)$ has property $(P)$.

\begin{definition} \label{def:ram} Let $A$ be a tree and $G < \aut(A)$.
\begin{itemize}
\item[(1)] A vertex incident to at least three edges is called a \emph{ramification point} \cite[2.1]{tits}.
\item[(2)] $G^+$ is the subgroup of $\aut(A)$ generated by the pointwise stabilizers in $G$ of edges \cite[4.5]{tits}: \[G^+ = \left\langle \stab^G(x,y) : \{x,y\}\in\Are(A) \right\rangle.\] 
\end{itemize}
\end{definition}

\begin{lemma} \label{lem:equ}
Every element of $G^+$ is either a rotation or a translation with an even translation length.
\end{lemma}
\begin{proof}
Consider the equivalence relation $E$ on $\Som(A)$: \[E(x,y)\ \ \Longleftrightarrow\ \ \text{the distance from $x$ to $y$ is even}.\] Every pointwise stabilizer of an edge fixes each $E$-class setwise, so $G^+$ preserves each $E$-class setwise. On the other hand, only rotations and translations with even translation lengths preserve each $E$-class setwise.
\end{proof}


For subsets $A,B$ of a group $G$, by $A^B$ we denote the set $\left\{a^b : a\in A, b\in B\right\}$, where $a^b=b^{-1}ab$, and for $n\in\N$, by $A^n$ we denote $\underbrace{A\cdot \ldots \cdot A}_{n\text{ times}}$.

\begin{definition}{\cite[2.3]{tits}}
Suppose $A$ is a tree. A subtree $A'$ of $A$ is called a \emph{half-tree} if $A'$ is a connected component of the forest obtained from $A$ by removing an edge (notice that, the obtained space has two connected components).
\end{definition}

\begin{proposition} \label{prop:tits}
Let $A$ be a tree, which is not a two-way infinite chain, and let $G<\aut(A)$. Assume that $G$ has property $(P)$ and that $G$ does not leave invariant any nonempty proper subtree of $A$, or any end of $A$.
\begin{enumerate}
\item For every nontrivial rotation $g\in G^+$, the set $g\cdot g^{G^+}$ contains a translation.
\item For every translation $g\in G^+$, the set $g^{G^+}\cdot \left(g^{-1}\right)^{G^+}$ contains the pointwise stabilizers of all half-trees, that is the set $\bigcup_{A' \subset A \text{ half-tree }} \stab^G(A')$.
\end{enumerate}
\end{proposition}
\begin{proof} The proof is a modification of the proof of \cite[Theorem 4.5]{tits}. 

$(1)$ Fix a nontrivial rotation $g\in G^+$. We use the following fact \cite[Lemma 4.4]{tits}: if $X$, $Y$ are nontrivial subgroups of $\aut(A)$, and $X$ normalizes $Y$, and $X$ does not leave invariant any nonempty proper subtree of $A$ or any end of $A$, then the same is true for $Y$; that is, $Y$ does not leave invariant any nonempty proper subtree of $A$ or any end of $A$. By applying this fact twice to $G$, to $G^+$ and to $\langle g^{G^+} \rangle$ we obtain that $\langle g^{G^+} \rangle$ does not leave invariant any nonempty proper subtree of $A$ or any end of $A$. Therefore, by Lemma \ref{lem:stab} applied to $X=g^{G^+}$, we find that the set $g^{G^+}\cdot g^{G^+} = \left( g\cdot g^{G^+} \right)^{G^+}$ contains a translation, so also $g\cdot g^{G^+}$ contains a translation.

$(2)$ Fix a translation $g\in G$. Let $A'$ and $A''$ be two half-trees obtained by removing some edge from $A$.
%
%
\begin{claim}
There is $h\in G^+$ such that the axis of $g^h$ is contained in $A'$.
\end{claim}
\begin{proof}[Proof of the claim]
As $G^+$ does not leave invariant any proper nonempty subtree or any end of $A$, the existence of $h$ follows by a part of the proof of \cite[Theorem 4.5]{tits}. We give the details.

Suppose $D$ is the axis of $g$. Without loss of generality we may assume that \[D \cap A'\neq \emptyset.\] Indeed, take an arbitrary vertex $s\in \Som(D)$. By \cite[Lemma 4.1]{tits}, there is $h\in G^+$ with $h^{-1}(s)\in \Som(A')$. Then $g^h$ has as axis $h^{-1}[D]$, and $h^{-1}[D]\cap A'\neq\emptyset$.

We may also assume that $D\not\subseteq A'$.

It is enough to find $h\in G^+$ such that $h[D]\subseteq A'$ (then the axis of $g^{h^{-1}}$ is $h[D]$).

Let $b'$ and $b''$ be the two ends of $A$ induced by $D$, such that $b'$ is an end of $A'$ and $b''$ is an end of $A''$. Since $G^+$ does not stabilize $b'$, $b''$ and $\{b',b''\}$ (otherwise $G^+$ leaves $D$ invariant and $A\ne D$), there is $g'\in G^+$ with \[g'(b'')\not\in\{b',b''\};\] indeed if $h(b'')=b'$ for some $h\in G^+$, then take $g''\in G^+$ with $g''(b')\not\in\{b',b''\}$ and put $g'=g''\circ h$. Denote by $\pi \colon \Som(A) \to \Som(D)$ a projection from $A$ to $D$, so $\dist(x,D) = \dist(x,\pi(x))$ for each $x\in A$. Since $b''\not\in\left\{{g'}^{-1}(b'),{g'}^{-1}(b'')\right\}$, the projection under $\pi$ of ${g'}^{-1}[D]$ has a finite number of vertices after intersecting with $A''$. Thus, there exists $n\in\Z$ such that $g^n\left[\pi\left[{g'}^{-1}[D]\right]\right]$ is included in $A'$. Then $g^n\left[g'^{-1}[D]\right]\subseteq A'$.
\end{proof}

Let $D\subseteq A'$ be the axis of $g^h$ from the claim. By \cite[Lemma 4.3]{tits}, using the assumption that $G$ has property $(P)$, \[\stab^G(D) = \left\{g^hf\left(g^h\right)^{-1}f^{-1} : f\in\stab^G(D)\right\}.\] Hence \[\stab^G\left(A'\right) < \stab^G(D) = g^h\cdot {\left(g^h\right)^{-1}}^{\stab^G(D)} \subseteq g^{G^+}\cdot \left(g^{-1}\right)^{G^+}.\]
\end{proof}

We recall from \cite[Section 5]{tits} a convenient way to describe trees. Let $I$ be a set of ``colors'' and \[f\colon \Som(A) \to I\] a coloring function. Define a group of automorphisms preserving $f$ as \[\autf(A) = \{\alpha \in \aut(A) : f \circ \alpha = f\}.\] We say that $f$ is \emph{normal} if $f$ is onto and for every $i\in I$, $\autf(A)$ is transitive on $f^{-1}[i]$. Clearly, for every coloring function $f$ there is a normal coloring function $f'$, possibly with a different set of colors, such that $\autf(A) = \aut_{\it f'}(A)$.  Hence we may always assume that $f$ is normal. 

It is easy to see that $\autf(A)$ has the property $(P)$.

Let $(A, f\colon \Som(A) \to I)$ be an arbitrary colored tree, with $f$ normal. Define a function \[a\colon I\times I \to \card\] as follows: take an arbitrary $x\in f^{-1}[i]$ and set \[a(i,j) = \left|\left\{y\in f^{-1}[j] : \{x,y\} \in \Are(A)\right\}\right|.\] Since $f$ is normal, the value $a(i,j)$ does not depend on the choice of $x$ from $f^{-1}[i]$. Functions $a$ arising in this way can be characterized by two conditions \cite[Proposition 5.3]{tits}:
\begin{itemize}
\item[(1)] if $a(i,j) = 0$, then $a(j,i) = 0$
\item[(2)] the graph $G(a) = (I, E)$, where $E = \{\{i,j\} \subseteq I : a(i,j)\neq 0\}$, is connected.
\end{itemize}
If a function $a\colon I\times I \to \card$ has properties $(1)$ and $(2)$, then there is a colored tree $A$ with a normal coloring function $f$ such that for every $x\in f^{-1}[i]$, \[a(i,j) = \left|\left\{ y \in f^{-1}[j] : \{x,y\} \in \Are(A) \right\}\right|.\] We say, then, that $a$ is a \emph{code} of the colored tree $A$. We note also \cite[5.7]{tits} that if $1\not\in a[I\times I]$, then $\autf(A)$ does not leave invariant any nonempty proper subtree or any end; hence by \cite[Theorem 4.5]{tits} $\autf^+(A)$ is a simple group.

An element $i\in I$ is a \emph{ramification color}, if $i=f(x)$, for some ramification point $x\in\Som(A)$. The set of all ramification colors is denoted by $I^{\ram}$.

By the \emph{color} of a chain (possibly infinite) we mean the sequence of colors of its vertices.

We will use the following fact from \cite{tits}.

\begin{proposition}{\cite[6.1]{tits}} \label{prop:ram} Let $A$ be a colored tree. The stabilizers of all ramification points are contained in $\autf^+(A)$, and $\autf^+(A)$ is generated by them, that is \[\autf^+(A) = \left\langle \stab^{\autf(A)}(r) : r\in\Som(A) \text{ is a ramification point} \right\rangle.\]
\end{proposition}

Let $I^+$ denote the set of orbits of $\autf^+(A)$ on $\Som(A)$. We write $f^+ \colon \Som(A) \to I^+$ for the induced quotient map.

\begin{proposition} \label{prop:quo}  
\begin{itemize}
\item[(1)] $f^+$ is normal and $f^+$ refines $f$, i.e. if $f^+(x) = f^+(y)$, then $f(x)=f(y)$.
\item[(2)] $\autf^+(A) = \aut_{f^+}(A)$
\end{itemize}
\end{proposition}
\begin{proof}
$(1)$ and the inclusion $\subseteq$ in $(2)$ are obvious. If $\alpha\in\aut_{f^+}(A)$ and $r\in\Som(A)$ is a ramification point, then $\alpha(r) \in\autf^+(A)\cdot r$. Thus, by Proposition \ref{prop:ram}, $\alpha\in \stab^{\autf(A)}(r) \cdot \autf^+(A)=\autf^+(A)$.
\end{proof}

\section{Bounded simplicity of $\autf^+(A)$}

We begin with the criterion for bounded simplicity of a group acting on a tree.




\begin{lemma} \label{lem:rami}
Assume that $(A, f\colon \Som(A) \to I)$ is a colored tree, $f$ is normal and the group $\autf^+(A)$ is nontrivial. 
\begin{enumerate}
\item Every nontrivial rotation from $\autf^+(A)$ fixes some ramification point and is a composition of two elements from $\bigcup_{\{x,y\}\in\Are(A)} \stab^{\autf(A)}(x,y)$.
\item Suppose that $G<\autf^+(A)$, $G^+$ is nontrivial, and that $G$ has property $(P)$ and does not leave invariant any nonempty proper subtree of $A$ or any end of $A$. Then $G^+$ is boundedly simple if and only if there is $N\in\N$ such that every translation from $G^+$ is the product of $N$ elements from $G^+$ each of which fixes poitwise a half-tree; in such case $G^+$ is $4N$-boundedly simple.
\end{enumerate}
\end{lemma}
\begin{proof} Note that, since $\autf^+(A)$ is nontrivial, $A$ is not a two-way infinite chain.

$(1)$ By \cite[6.1]{tits}, if $\alpha\in\autf^+(A)$ stabilizes a ramification point, then $\alpha$ is a product of two elements from $\bigcup_{\{x,y\}\in\Are(A)} \stab^{\autf(A)}(x,y)$. We prove that every rotation $\alpha\in\autf^+(A)$ fixes a ramification point. Introduce the following equivalence relation $E$ on the set of all ramification points of $A$: for ramification points $r_1\ne r_2$, $E(r_1,r_2)$ holds if and only if on the chain joining $r_1$ and $r_2$ there is an odd number of ramification points. $E$ has exactly two equivalence classes. Each of the equivalence classes of $E$ is invariant under $\autf^+(A)$. Therefore, a rotation which does not preserve any ramification point is not in $\autf^+(A)$.
%

$(2)$ $\Rightarrow$ is clear, since there exists at least one nontrivial element from the pointwise stabilizer of some half-tree: $G^+$ is nontrivial and $G$ has property $(P)$, so for any edge $\{x,y\}$ there are half-trees $A'$, $A''$ such that $\stab^G(x,y) = \stab^G(A')\cdot\stab^G\left(A''\right)$.

$\Leftarrow$ First, let $g\in G^+$ be an arbitrary translation. Then by Proposition \ref{prop:tits}$(2)$, the set $\left(g^{G^+}\cdot \left( g^{-1} \right)^{G^+}\right)^N$ contains all translations from $G^+$. By $(1)$, every rotation from $G^+$ is a product of two elements each fixing two edges. Thus, by property $(P)$, it is a product of four elements each fixing pointwise half-trees, so is in $\left(g^{G^+}\cdot \left( g^{-1} \right)^{G^+}\right)^4$ (again by \ref{prop:tits}$(2)$). Now, if $g\in G^+$ is a nontrivial rotation, then by Proposition \ref{prop:tits}$(1)$, $gg^h$ is a translation, for some $h\in G^+$, and in this case $G^+ = \left(\left(gg^h\right)^{G^+}\cdot {\left( gg^h \right)^{-1}}^{G^+}\right)^N$.
\end{proof}

We now define the main ingredient of the later proofs in this paper, that is the notion of the \emph{type of a translation}. We associate with each translation, a finite sequence of colors.

\begin{definition}
Let $(A, f\colon \Som(A) \to I)$ be an arbitrary colored tree and $\alpha\in\autf(A)$ be a translation along a two-way infinite chain $C$. Take an arbitrary vertex $x\in \Som(C)$ and a subchain $(x_1,\ldots,x_{n+1})$ of $C$ such that $x = x_1,\ldots,x_{n+1} = \alpha(x)$. Define $i_1:=f(x_1),\ldots,i_n:=f(x_n)$, noting that $f(x_{n+1})=i_1$. Then we say that the set \[[i_1,\ldots,i_n] = \left\{(i_1,\ldots,i_n), (i_2,\ldots,i_n,i_1),\ldots,(i_n,i_1,\ldots,i_{n-1})\right\}\] of all cyclic shifts of the sequence $(i_1,\ldots,i_n)$ is the \emph{type} of the translation $\alpha$.
\end{definition}

Any two translations which are conjugate have the same type. We calculate types of some translations: a composition of two rotations and a composition of a rotation and a translation.

\begin{figure}
\centering
\includegraphics{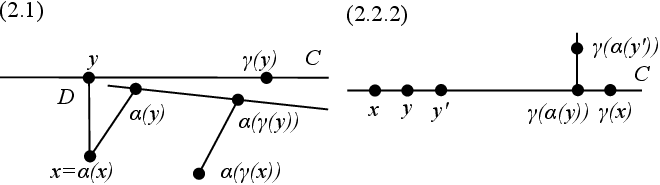}
\caption{Composition of translation and rotation}
\label{trans.rot}
\end{figure}

\begin{lemma} \label{lem:an}
Let $\alpha,\beta\in\autf(A)$ be rotations such that $\fix(\alpha)\cap\fix(\beta)=\emptyset$ and let $\gamma\in\autf(A)$ be a translation.
\begin{itemize}
\item[(1)] Assume that $\alpha(x)=x$, $\beta(y)=y$ for some $x,y\in\Som(A)$, and on the chain $D$ from $y$ to $x$ the only vertex fixed by $\alpha$ is $x$ and the only vertex on $D$ fixed by $\beta$ is $y$. If the color of $D$ is $(i_1,\ldots,i_n)$, where $n\geq 2, f(y)=i_1$, and $f(x)=i_n$, then the type of $\alpha\circ\beta$ is \[\tau_2 =[i_1,i_2,\ldots,i_{n-1},i_n,i_{n-1},\ldots,i_2].\] Moreover, every translation of the type $\tau_2$ is a composition of two rotations.

\item[(2)] Assume that $\gamma$ is a translation of the type $\tau=[i_1,j_2,\ldots,j_m]$, where $m\geq 2$, along a two-way infinite chain $C$ and $x$ is the vertex fixed by $\alpha$ which is closest to $C$.
	\begin{itemize}
	\item[(2.1)] Assume that $x$ lies outside $C$ (see Figure \ref{trans.rot}). Let $y\in C$ be the closest vertex to $x$ in $C$, that is a projection of $x$ on $C$. Let $D$ be the chain from $y$ to $x$. Let the color of $D$ be $(i_1,\ldots,i_n)$ with $f(y)=i_1$ and $f(x)=i_n$.
Then $\alpha\circ\gamma$ and $\gamma\circ\alpha$ are translations of the type \[\tau_3 = [i_1,i_2,\ldots,i_{n-1},i_n,i_{n-1},\ldots,i_2,i_1,j_2,\ldots,j_m].\] Also, every translation of the type $\tau_3$ is a composition of a rotation, and a translation of the type $\tau$. 
	\item[(2.2)] Assume that $x$ lies on $C$. Let $D$ be the chain from $x$ to $\gamma(x)$ and assume that $f(x)=i_1$. Let $y$ be a vertex from $D$ adjacent to $x$ (so $f(y)=j_2$).
\begin{itemize}
\item[(2.2.1)] If $\gamma(\alpha(y))$ lies outside $D$, then $\gamma\circ\alpha$ is a translation of the same type as $\gamma$, that is of the type $\tau$.
\item[(2.2.2)] Assume that $\gamma(\alpha(y))$ is on $D$ (so $j_2 = j_m$). Let $y'\neq x$ be a vertex from $D$, adjacent to $y$ (so $f(y')=j_3$). If $\gamma(\alpha(y'))$ is outside $D$, then $\gamma\circ\alpha$ is a translation of the type $[j_2,\ldots,j_{m-1}]$.
\end{itemize}
\end{itemize}
For $\gamma\circ\alpha$ exactly one of the following statements is true.
\begin{enumerate}
\item $\gamma\circ\alpha$ is a translation of the type being the subtype of $\tau$.
\item $\gamma\circ\alpha$ is a rotation. In this case $m$ is even and $j_2 = j_{m},\ j_3 = j_{m-1},\ \ldots,\ j_{\frac{m}{2} - 1} = j_{\frac{m}{2} + 2}$. Thus $\gamma\circ\alpha$ stabilizes vertex of type $j_{\frac{m}{2} + 1}$.
\item $\gamma\circ\alpha$ is an inversion. In this case $m$ is odd.
\end{enumerate}
Since $\alpha\circ\gamma = (\gamma\circ\alpha)^{\alpha^{-1}}$, the same applies to $\alpha\circ\gamma$.
\end{itemize}
\end{lemma}
\begin{proof} By applying $(\spadesuit)$ from Lemma \ref{lem:trans} to: 
\begin{itemize}
\item $y$, $x$, $\alpha(\beta(y))$, $\alpha(\beta(x))$, in $(1)$,
\item $y$, $x$, $\alpha(\gamma(y))$, $\alpha(\gamma(x))$ (see Figure \ref{trans.rot}), in $(2.1)$,
\item $x$, $y$, $\gamma(x)=\gamma(\alpha(x))$ and $\gamma(\alpha(y))$, in $(2.2.1)$,
\item $y$, $y'$, $\gamma(\alpha(y))$ and $\gamma(\alpha(y'))$ (see Figure \ref{trans.rot}), in $(2.2.2)$,
\end{itemize}
we have in $(1)$ that the type of $\alpha\circ\beta$ is $\tau_2$, and in $(2)$ that either the type of $\gamma\circ\alpha$ is the subtype of $\tau$, or $\gamma\circ\alpha$ is a rotation.

We prove that every translation of type $\tau_2$ is a composition of two rotations. The case of a composition of a translation and a rotation is similar. Suppose $a$ is a code of $A$. Let $\delta$ be a translation of type $\tau_2$ along the chain $C$. Then $a(i_1,i_2)$ and $a(i_n,i_{n-1})$ are at least 2. Hence there are rotations $\alpha,\beta\in\autf(A)$ such that $\alpha\circ\beta$ has the same type and translation chain as $\delta$. Then $\delta' = \delta\circ(\alpha\circ\beta)^{-1}$ fixes $C$, so $\delta = (\delta'\circ\alpha)\circ\beta$ is a composition of two rotations.
\end{proof}



An $(n,m)$-regular or biregular tree, denoted by $A_{n,m}$, is a 2-colored tree with the following code: \[a(0,0) = a(1,1) = 0,\ \ a(0,1) = n,\ \ a(1,0) = m,\] where $I=\{0,1\}$ and $n$, $m$ are some cardinal numbers $\geq3$. Intuitively, in a biregular tree every vertex is black or white, every white vertex is connected with $n$ black vertices and every black vertex is connected with $m$ white vertices (if we allow $n=2$ and $m\geq3$, then after removing vertices of color $0$ we get the subdivision of an $m$-regular tree $A_{m,m}$).

\begin{theorem} \label{thm:bireg}
Suppose that $n, m\geq 3$ are cardinals. The group $\aut^+(A_{n,m})$ is $8$-boundedly simple. Moreover, if $m=n$, then $\left[\aut(A_{n,m}):\aut^+(A_{n,m})\right]= 2$; if $m\ne n$, then $\aut(A_{n,m})^+=\aut(A_{n,m})$.
\end{theorem}
\begin{proof} The moreover part is obvious. We prove the first part. Clearly, $\aut^+(A_{n,m})$ has property $(P)$ and $1\not\in a[I\times I]$, so $\aut^+(A_{n,m})$ leaves no nonempty proper subtree of $A_{n,m}$ invariant and does not stabilize any end of $A_{n,m}$. Also, $\aut^+(A_{n,m})$ consists of all translations with even translation lengths and all rotations. Thus, $\aut^+(A_{n,m})$ has exactly two orbits on the set $\Som(A)$, and every element of $\aut(A_{n,m})$ that leaves these two orbits invariant is in $\aut^+(A_{n,m})$. Therefore any two translations of the same translation length are conjugate by an element of $\aut^+(A_{n,m})$. Also, as for every $n$ there exists a product of two elements stabilising half-trees, which is a translation of the length $2n$. Hence, every translation from $\aut^+(A_{n,m})$ is the product of two elements stabilising half-trees. Now, the conclusion follows by Lemma \ref{lem:rami}.
\end{proof}

\begin{definition} \label{def:almost}
A \emph{uniform subdivision of $(n,m)$-regular tree} is the $(n,m)$-regular tree subdivided (in an equivariant way) by non-ramification points. Namely, it is the tree with the set of colors $I = \{0,\ldots,k\}$ and the following code:  $a(0,1)=n$, $a(k,k-1)=m$ and $a(i,i+1) = a(i,i-1) = 1$ for $i\in I \setminus \{0,k\}$. For all other pairs $(p,q)$ from $I^2$, $a$ has value $0$.
\end{definition}

If $A$ is a uniform subdivision of an $(n,m)$-regular tree $A_{n,m}$, then $\autf^+(A) \cong \aut^+(A_{n,m})$. Hence $\autf^+(A)$ is $8$-boundedly simple.

Apart from $A_{n,m}$, there are no other colored trees $A$ with boundedly simple groups $\autf^+(A)$ and with the property that $\autf^+(A)$ leaves no nonempty proper subtree of $A$ invariant (Theorem \ref{thm:main}). Proposition \ref{prop:autf} is the main technical step in the proof of this fact. We prove that, if $\autf^+(A)$ is boundedly simple, then some particular configuration in the code of $A$ is forbidden. In the proof of \ref{prop:autf} we use some combinatorial argument, describing the complexity of distances of colors in types.

\begin{definition} \label{def:iseq}
For $i\in I$ and type $t = [i_1,\ldots,i_n]$ define the \emph{$i$-sequence of $t$} in the following way.
\begin{itemize}
\item If there is no occurrence of $i$ in $t$, then the $i$-sequence of $t$ is empty.
\item Let $i_k$ be the first occurrence of $i$ in $(i_1,\ldots,i_n)$. The $i$-sequence of $t$ is a sequence (modulo all cyclic shifts) of distances between consecutive occurrences of $i$ in the sequence $(i_k,i_{k+1},\ldots,i_{n-1},i_n,i_1,\ldots,i_k)$.
\end{itemize}
\end{definition}

\begin{definition}
For any $i\in I$ and any type $t$ define $O(t,i)$ as the number of integers that appear an odd number of times in the $i$-sequence of $t$.
\end{definition}

\begin{lemma} \label{lem:occ}
Let $t$ be the type of a translation which is a composition of $K$ rotations. Then $O(t,i) \leq 4K - 6$.
\end{lemma}

\begin{proof}
We prove the lemma by induction on $K$. 

Let $K=2$. By Lemma \ref{lem:an}$(1)$, $t$ is of the form $[i_1,i_2,\ldots,i_{n-1},i_n,i_{n-1},\ldots,i_2]$. In the Table \ref{table:two} we describe all possibilities for the shape of the $i$-sequence of $t$.
\begin{table}
\caption{$i$-sequence of the composition of two rotations}
\label{table:two}
\centering
\begin{tabular}{| c | c |}
\hline
 Case                            & $i$-sequence\\ \hline
 $i_1 = i$ and $i_n = i$         &  $\left[m_1,m_2,\ldots,m_{\frac{N}{2}},m_{\frac{N}{2}},\ldots,m_2,m_1\right]$ \\ \hline
 $i_1 = i$ and $i_n \neq i$      &  $\left[m_1,m_2,\ldots,m_{\frac{N-1}{2}},m_{\frac{N+1}{2}},m_{\frac{N-1}{2}},\ldots,m_2,m_1\right]$\\ \hline
 $i_1 \neq i$ and $i_n = i$      &  $\left[m_1,m_2,\ldots,m_{\frac{N-1}{2}},m_{\frac{N-1}{2}},\ldots,m_2,m_1,m_0\right]$\\ \hline
 $i_1 \neq i$ and $i_n \neq i$   &  $\left[m_1,m_2,\ldots,m_{\frac{N-2}{2}},m_{\frac{N}{2}},m_{\frac{N-2}{2}},\ldots,m_2,m_1,m_0\right]$\\
\hline
\end{tabular}
\end{table}
In all cases $O(t,i) \leq 2$.

Let $t$ be the type of the composition $\tau = \tau_1\circ\ldots\circ\tau_{K+1}$ of $K+1$ rotations. Put $\rho = \tau_1\circ\ldots\circ\tau_{K}$. If $\rho$ is a rotation, then $\tau = \rho \circ \tau_{K+1}$ is the composition of two rotations and we may use the induction hypothesis. Otherwise, $\rho$ is a translation along some two-way infinite chain $C'$. Let $s$ be the type of $\rho$ and $x_{K+1}$ be the vertex fixed by $\tau_{K+1}$, which is nearest to $C'$. There are two cases: $x_{K+1}$ is in $C'$ or not.

Assume first that $x_{K+1} \not\in C'$, i.e. that the case $(2.1)$ from Lemma \ref{lem:an} holds. Then \[s = [i_1,j_2,\ldots,j_m]\ \text{ and }\ t = [i_1,i_2,\ldots,i_{n-1},i_n,i_{n-1},\ldots,i_2,i_1,j_2,\ldots,j_m],\] for some $n, m\geq 2$. Let $N_1$ and $N_2$ be the numbers of occurrences of $i$ in $(i_1,\ldots,i_n,\ldots,i_2)$ and $(i_1,j_2,\ldots,j_m)$ respectively. Let $[n_1,n_2,\ldots,n_{N_2}]$ denote the $i$-sequence of $s$. Again, there are four possibilities for the shape of the $i$-sequence of $t$ (presented in the Table \ref{table:many}).
\begin{table}
\caption{$i$-sequence of $t$}
\centering
\begin{tabular}{| c | c |}
\hline
 Case                            & $i$-sequence\\ \hline
 $i_1 = i$ and $i_n = i$         &  $\left[m_1,\ldots,m_{\frac{N_1}{2}},m_{\frac{N_1}{2}},\ldots,m_1,n_1,n_2,\ldots,n_{N_2}\right]$ \\ \hline
 $i_1 = i$ and $i_n \neq i$      &  $\left[m_1,\ldots,m_{\frac{N_1-1}{2}},m_{\frac{N_1+1}{2}},m_{\frac{N_1-1}{2}},\ldots,,m_1,n_1,n_2,\ldots,n_{N_2}\right]$\\ \hline
 $i_1 \neq i$ and $i_n = i$      &  $\left[m_1,\ldots,m_{\frac{N_1-1}{2}},m_{\frac{N_1-1}{2}},\ldots,m_1,m_0,n_1,n_2,\ldots,n_{N_2-1},n'_{N_2}\right]$\\ \hline
 $i_1 \neq i$ and $i_n \neq i$   &  $\left[m_1,\ldots,m_{\frac{N_1-2}{2}},m_{\frac{N_1}{2}},m_{\frac{N_1-2}{2}},\ldots,m_1,m_0,n_1,n_2,\ldots,n_{N_2-1},n'_{N_2}\right]$\\
\hline
\end{tabular}
\label{table:many}
\end{table}
By the induction hypothesis $O(s,i) \leq 4K - 6$. Therefore, by the definition of $O(t,i)$, in the worst (i.e. fourth) case we have \[O(t,i) \leq O(s,i) + 4 \leq 4(K+1) - 6.\]

Assume now that $x_{K+1} \in C'$, i.e. the case $(2.2)$ from Lemma \ref{lem:an} holds. Then \[s = [i_1,i_2,\ldots,i_{n-1},i_n,i_{n-1},\ldots,i_2,i_1,j_2,\ldots,j_m]\ \text{ and }t = [i_1,j_2,\ldots,j_m]\ n, m\geq 2.\] Let $N_1$ and $N_2$ be the numbers of occurrences of $i$ in $(i_1,i_2,\ldots,i_{n-1},i_n,i_{n-1},\ldots,i_2)$ and $t$ respectively. We may assume that the $i$-sequence of $s$ is given by the Table \ref{table:many} (where $[n_1,n_2,\ldots,n_{N_2}]$ is the $i$-sequence of $t$). By the induction hypothesis $O(s,i) \leq 4K - 6$. We have to show that $O(t,i) \leq 4(K+1) - 6$. In the first case (i.e. $i_1 = i$ and $i_n = i$) $O(t,i) = O(s,i)$. In the second case $O(t,i) \leq O(s,i) + 1$. In the third case $O(t,i) \leq O(s,i) + 3$. In the fourth case $O(t,i) \leq O(s,i) +4 \leq 4(K+1) - 6$.
\end{proof}

\begin{figure}
\centering
\includegraphics{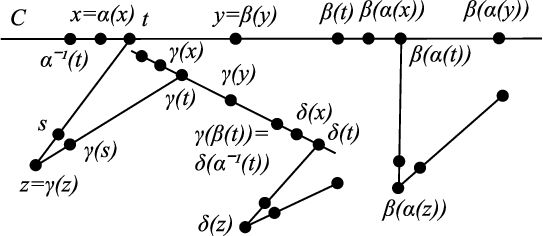}
\caption{Composition of three rotations}
\label{3.rot}
\end{figure}

\begin{proposition} \label{prop:autf}
Assume that $A$ is a colored tree and $\autf^+(A)$ is nontrivial and boundedly simple. Let $\alpha$, $\beta\in \autf^+(A)$ be rotations such that $\fix(\alpha)\cap\fix(\beta)=\emptyset$. Suppose that for three different ramification points $x,y,z\in \Som(A)$,
\begin{itemize}
\item $\alpha(x) = x$, $\beta(y) = y$ and on the chain from $x$ to $y$ the only vertex fixed by $\alpha$ is $x$ and the only vertex fixed by $\beta$ is $y$; i.e. $\beta\circ\alpha$ is a translation along a two-way infinite chain $C$ (see Figure \ref{3.rot}),
\item $t$ is the projection of $z$ onto $C$ and $s$ is a vertex adjacent to $z$ lying on the chain from $z$ to $t$,
\item on the chains from $x$ to $y$ and from $s$ to $t$ there are no vertices of color $f(z)$ (so also on $C$ there are no such vertices).
\end{itemize}
Let $\gamma \in \autf^+(A)$ with $\gamma(z)=z$, then $\gamma(s) = s$.
\end{proposition}

\begin{proof} By Proposition \ref{prop:quo}, there is a normal function $f^+ \colon \Som(A) \to I^+$ with $\autf^+(A) = \aut_{f^+}(A)$. We may assume further that $f = f^+$ and $I = I^+$.

Suppose, contrary to our claim, that $\gamma(s)\neq s$. For each $K\in\N$ we construct a composition of some rotations which cannot be written as a composition of $K$ rotations. Then, Lemma \ref{lem:rami} implies that $\autf^+(A)$ is not boundedly simple.

Our situation is described by Figure \ref{3.rot}. We may assume that $t$ belongs to the chain in $C$ from $x$ to $y$ (if $t$ belongs to the chain from $(\beta\circ\alpha)^n(x)$ to $(\beta\circ\alpha)^n(y)$, for some integer $n$, then just take $z:=(\beta\circ\alpha)^{-n}(z)$ and the conjugate $\gamma := \gamma^{(\beta\circ\alpha)^n}$).

Denote by $u$, $v$ and $w$ sequences of colors, corresponding to chains in Figure \ref{3.rot2}. Namely, let
\begin{itemize}
\item $u$ corresponds to the chain from $t$ to $\gamma(t)$ (through $z$) without the last term of color $f(t)$,
\item $v$ --- from $\alpha^{-1}(t)$ to $t$ (through $x$) without the last term of color $f(t)$,
\item $w$ --- from $t$ to $\beta(t)$ (through $y$) without the last term of color $f(t)$.
\end{itemize}

Note that $v$ or $w$ might be empty, but the chains $u$ and $vw$ (the concatenation of $v$ and $w$) are always nonempty.

By Lemma \ref{lem:an}, sequences $u$, $v$ and $w$ have the form \[(c_1,c_2,\ldots,c_{r-1},c_r,c_{r-1},\ldots,c_2), \tag{$\blacklozenge$}\] where $c_1=f(t)$.

For example, by Lemma \ref{lem:an}$(2.1)$, translations $\alpha\circ\beta$, $\beta\circ\alpha$ have type $[v,w]$, translations $\delta = \gamma\circ\beta\circ\alpha$, $\beta\circ\alpha\circ\gamma$, $\gamma\circ\alpha\circ\beta$, $\alpha\circ\beta\circ\gamma$ have type $[v,u,w]$ and the chain $C$ has color $(\ldots vwvw \ldots)$ (see Figure \ref{3.rot2}).

\begin{figure}
\centering
\includegraphics{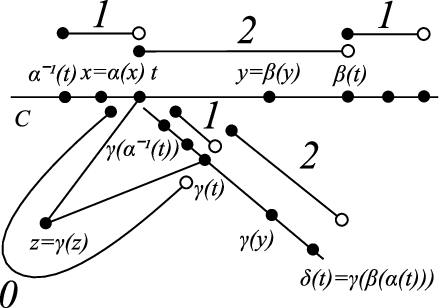}
\caption{Types of chains}
\label{3.rot2}
\end{figure}

Define by induction the following sequences 
\begin{itemize}
\item $t_2 = (v,w,v,u,(v,w)^2,v,u)$,
\item $t_{n+1} = (v,w,v,u,(v,w)^{2n-1},\ t_n,(v,w)^{2n-1},v,u)$, for $n\geq 2$. 
\end{itemize}

Note that each $t_n$ induces a two-way infinite chain in $A$ of the color $(\ldots t_nt_nt_n \ldots)$. In fact, vertex $t$ is a ramification point as a joining point of chains $u$, $v$ and $w$. Hence by $(\blacklozenge)$, for each automorphism $\rho\in\aut(A)$, one can find in $A$ segments starting at $\rho(t)$ of type $u$, $v$ and $w$.

Let $\alpha_n$ be a translation of type $[t_n]$ ($\alpha_n$ is the translation along the chain $(\ldots t_nt_nt_n \ldots)$). Then $\alpha_n$ is the composition of $n$ rotations from $\autf^+(A)$. In particular, by Lemma \ref{lem:an}$(1)$, $[t_2] = [w,v,u,(v,w)^2,v,u,v]$ is a type of composition of two rotations from $\autf^+(A)$ (because $x$ and $y$ are ramification points). Also $\alpha_{n+1}$ has the type \[[v,w,v,u,(v,w)^{2n-1},\ t_n,(v,w)^{2n-1},v,u] = [(v,w)^{2n-1},v,u,v,w,v,u,(v,w)^{2n-1},t_n],\] being (by Lemma \ref{lem:an}$(2.1)$) the type of the composition of a translation of type $t_n$ and a rotation. Hence, $\alpha_{n+1}$ is a composition of $n+1$ rotations.

The proof will be completed by showing that $\alpha_n$ cannot be written as a composition of less than $\frac{n}{2}+1$ rotations.

We compute the $f(z)$-sequence of $t_n$ (see Definition \ref{def:iseq}). Note that, by the assumption, $f(z)$ appears once only in the chain $u$. However, first we compute the $u$-sequence for $t_n$ (regarding $u$ as an additional color). The $u$-sequence for $[t_2]$ is $[6,4]$ and for $[t_3]$, it is $[10,6,8,4]$. It can be proved by induction that the $u$-sequence for $[t_{n+1}]$ is \[(4n+2,4n-2,\ldots,14,10,6,8,12,\ldots,4n-4,4n,4). \tag{$\bigstar$}\] 

Let $p$ be the length of the chain $uv$ and $q$ the length of the chain $vw$ ($p$, $q$ are even and at least $2$). Note that the $u$-sequence for $t' = [u,(v,w)^n,1]$ is $(2n+2)$ and the $f(z)$-sequence for $t'$ is $(p + nq)$. Therefore by $(\bigstar)$, the $f(z)$-sequence for $[t_{n+1}]$ is \[[p+(2n)q,p+(2n-2)q,\ldots,p+4q,p+2q, p+3q,\ldots,p+(2n-3)q,p+(2n-1)q,p+q].\]

The $f(z)$-sequence for $t_n$ has no multiple occurrences of any value and $t_n$ has $2n-2$ occurrences of $f(z)$. In fact $O(t_n,f(z)) = 2n-2$. If $t_n$ is the type of the composition of $K$ rotations, then by Lemma \ref{lem:occ}, $2n-2 \leq 4K - 6$, so $K\geq \frac{n}{2}+1$. This finishes the proof of Proposition \ref{prop:autf}.
\end{proof}

Proposition \ref{prop:autf} implies that for many trees $A$, the groups $\autf^+(A)$ are not boundedly simple. That is, after adding to ``an almost arbitrary'' tree $A$ one new color $k$, such that for some old color $j$, $a(k,j)\geq2$, we obtain a tree $A'$ where $\autf^+(A')$ is not boundedly simple.

\begin{corollary}
Assume that $(A, f\colon \Som(A) \to I)$ is a colored tree,  $f$ is normal and $\autf^+(A)$ does not stabilize any vertex and does not centralize any end. Extend the code $a$ of $A$ by adding one new color $I' = I \cup \{k\}$ ($k\not\in I$) to get a code $a'\supset a$ such that: $k$ is a ramification color, for every $i\in I$, $a'(i,k)=0$ if and only if $a'(k,i)=0$, and for some $j\in I$, $a'(k,j)\geq 2$. If $(A', f'\colon \Som(A') \to I')$ is a tree corresponding to $a'$, then $\aut_{f'}^+(A')$ is not boundedly simple.
\end{corollary}
\begin{proof} The tree $A'$ contains the subtree $A$ corresponding to $a$. Let $z$ be a vertex in $A'$ of color $k$ and let $s$ be a vertex in $A$ of color $j$ adjacent to $z$. Since $\autf^+(A)$ does not stabilize any vertex or any end, by \cite[Proposition 3.4]{tits} there exists a translation in $\autf^+(A)$ along a two-way infinite chain $C$ in $A$ and which is a composition of two rotations from $\autf^+(A)$. Let $t$ be the projection of $s$ onto $C$ in the tree $A$. Applying Proposition \ref{prop:autf} to $z$, $s$, $t$ and $C$, we conclude that $\aut_{f'}^+(A')$ is not boundedly simple (because there is $\gamma\in\aut_{f'}^+(A')$, such that $\gamma(z)=z$ and $\gamma(s)\neq s$).
\end{proof}

We characterize all colored trees $A$ with boundedly simple group $\autf^+(A)$ under the assumptions of \cite[Theorem 4.5]{tits}.


\begin{theorem} \label{thm:main}
Assume that $(A, f\colon \Som(A) \to I)$ is a colored tree and $\autf^+(A)$ is boundedly simple and nontrivial. Then $\autf^+(A)$ fixes some vertex of $A$, or leaves invariant some end of $A$, or leaves invariant a subtree $A' \subseteq A$, which is a uniform subdivision of an $(n,m)$-regular tree, for some $n,m\geq 3$.

In particular, if $\autf^+(A)$ leaves no nonempty proper subtree of $A$ invariant and stabilizes no end, then $A$ is a uniform subdivision of a biregular tree and $\autf^+(A)$ is $8$-boundedly simple.
\end{theorem}
\begin{proof} 
The second part of the theorem follows by the first part and Theorem \ref{thm:bireg}. We prove the first part.

By Proposition \ref{prop:quo}, there is a normal function $f^+ \colon \Som(A) \to I^+$ such that $\autf^+(A) = \aut_{f^+}(A)$. Let $a^+$ be the code for $(A, f^+\colon \Som(A) \to I^+)$. A short argument using Proposition \ref{prop:ram} yields that the procedure of removing non-ramification vertices does not change the group $\autf^+(A)$. Hence, without loss of generality, we may assume that all vertices of $A$ are ramification vertices. We may also assume that $\autf^+(A)$ is infinite, as by a Bruhat-Tits fixed point theorem, any finite group action on a tree has a global fixed point \cite[Theorem 15]{serre}.

We use the following property of translations from $\autf^+(A)$, due to Cong Chen.

\begin{lemma}[Cong Chen] \label{lem:cong}
Suppose $\alpha\in \autf^+(A)$ is a translation of type $t = [i_1,\ldots,i_n]$. For any two different colors $i,j\in I^+$, the number of occurrences of the sequence $(i,j)$ as a subsequence of consecutive terms in the sequence $(i_1,\ldots,i_n,i_1)$ equals the number of occurrences of $(j,i)$ in $(i_1,\ldots,i_n,i_1)$.
\end{lemma}
\begin{proof}
The proof goes along the lines of Lemma \ref{lem:an}.
\end{proof}

As a straightforward consequence of the lemma we have:

\begin{itemize}
\item[$(\clubsuit)$] Let $G(a^+)=(I^+, E)$ be the graph on $I^+$ where $E = \left\{\{i,j\} \subseteq I^+ : a^+(i,j)\neq 0\right\}$ (see Section 2). Then $G(a^+)$ has no cycles, i.e. $G(a^+)$ is a tree.
\end{itemize}

Indeed, suppose $c = (i_1,\ldots,i_n,i_1)$ is a cycle in $G(a^+)$ of pairwise distinct colors and $n\geq 3$. Then, there exists a chain $x_1,\ldots,x_n,x_{n+1}$ in $A$ of color $c$. Thus $x_1 = \alpha(x_{n+1})$ for some $\alpha\in \autf^+(A)$. However, $\alpha$ cannot be a rotation, since the colors $i_1,\ldots,i_n$ are pairwise distinct, and by the lemma, $\alpha$ cannot be a translation (because the type of $\alpha$ would be a subsequence of $(i_1,\ldots,i_n)$). This proves $(\clubsuit)$.

We may assume that $\autf^+(A)$ does not fix any vertex of $A$ or stabilize any end of $A$.

\begin{clm} \label{clm:0}
$\autf^+(A)$ does not leave invariant any edge of $A$.
\end{clm}
\begin{proof}
Suppose an edge $\{s,t\}$ is $\autf^+(A)$-invariant. The stabilizer $\stab^{\autf^+(A)}(s)$ has index at most 2 in $\autf^+(A)$, so is a normal subgroup. As $\autf^+(A)$ is simple, the vertex $s$ is $\autf^+(A)$-invariant, which is impossible.
\end{proof}

By Claim \ref{clm:0} and \cite[Proposition 3.4]{tits}, $\aut_{f^+}(A)$ contains some translation which is a composition of two rotations from $\aut_{f^+}(A)$ (by Lemma \ref{lem:trans}). Let $\alpha\in\aut_{f^+}(A)$ be a translation of the minimal possible translation length amongst all translations which are products of two rotations. 

Let \[[j_0,j_1,\ldots,j_{k-1},j_k,j_{k-1},\ldots,j_1],\ \ k\geq1\] be the type of $\alpha$ according to the coloring $f^+$. We may assume that $\alpha = \beta\circ\gamma$, for some $\beta, \gamma\in \aut_{f^+}(A)$ such that $\beta(x)=x$, $\gamma(y)=y$; also, that the colors of $x$ and $y$ are $j_0$ and $j_k$ respectively, and \[(x=x_0,x_1,\ldots,x_{k-1},x_k=y) \tag{\textasteriskcentered} \] is the chain in $A$ from $x$ to $y$ (so $f^+(x_i)=j_i$). Notice that $j_i\ne j_{i+1}$ for $0\leq i \leq k-1$, by Lemma \ref{lem:equ}. Define $n = a^+(j_0,j_1)$ and $m = a^+(j_k,j_{k-1})$. Clearly $n,m\geq 2$.

\begin{clm} \label{clm:1}
If $k\geq 2$, then $j_{i-1}\ne j_{i+1}$ for every $0<i<k$.
\end{clm}
\begin{proof}[Proof of Claim \ref{clm:1}]
Suppose that $j_{i-1}=j_{i+1}$. Then there exists $\beta'\in \aut_{f^+}(A)$ such that $\beta'(x_{i-1})=x_{i+1}$ and $\beta'(x_i)=x_i$ (it is enough to define inductively $\beta'$ as a permutation of $A$ preserving $f^+$). Hence by Lemma \ref{lem:an}$(1)$, the translation length of $\alpha' = \beta'\circ\gamma$ is smaller than the translation length of $\alpha$.
\end{proof}

\begin{clm} \label{clm:2}
If $k\geq 2$, then $a^+(j_1,j_2) = \ldots = a^+(j_{k-1},j_k) = 1$ and $a^+(j_{k-1},j_{k-2}) = \ldots = a^+(j_1,j_0) = 1$.
\end{clm}
\begin{proof}[Proof of Claim \ref{clm:2}]
Suppose that $a^+(j_i,j_{i+1})\geq 2$, for some $0<i<k$. Then instead of $\beta$ we may consider a nontrivial rotation $\beta'$ of $A$ fixing $x_i$ with $\beta'(x_{i+1})\ne x_{i+1}$. Note that by Proposition \ref{prop:ram}, $\beta'$ is in $\aut_{f^+}(A)$. We get the contradiction in the same way as in Claim \ref{clm:1}.
\end{proof}

\begin{clm} \label{clm:3}
The colors $j_0,j_1,\ldots,j_k$ are pairwise distinct.
\end{clm}
\begin{proof}[Proof of Claim \ref{clm:3}]
Fix $s, t\in \{1,\ldots,k-1\}$. Suppose that $j_s = j_t$. Then $\delta(x_s)=x_t$, for some $\delta\in\aut_{f^+}(A)$. If $\fix(\delta)\ne\emptyset$, then $\delta(x_i)=x_i$, for some $1< i < k-1$, which is impossible by Claim \ref{clm:1} (because then $\delta(x_{i-1})=x_{i+1}$). Therefore $\delta$ is a translation. Hence by Claim \ref{clm:2}, $\delta(x_{s+1})=x_{t+1}$, $\delta(x_{s-1})=x_{t-1}$, and more generally $\delta(x_n)=x_{n+t-s}$, for $0 \leq n \leq k - (t-s)$. Therefore $\delta(x_0)=x_{t-s}$, $\delta(x_{k-(t-s)})=x_k$ and $j_0 = j_{t-s}$. Let $\beta':= \beta^{\delta^{-1}}$. Then $\beta'(x_{t-s})=x_{t-s}$, so $\beta'\in \aut_{f^+}(A)$ (by Proposition \ref{prop:ram}). Thus $\alpha'=\beta'\circ\gamma$ has smaller translation length than $\alpha$, because $\beta'(x_k)\ne x_k$ (otherwise $\beta(\delta^{-1}(x_k))=\delta^{-1}(x_k)$, so $\beta(x_{k-(t-s)})=x_{k-(t-s)}$ which is impossible). The proofs of $j_s\ne j_0$, $j_s\ne j_k$ are similar.

Suppose $j_0=j_k$. Then for some $\delta\in\aut_{f^+}(A)$ we have $\delta(x_0)=x_k$. As in the previous case $\delta$ cannot be a rotation, so it must be a translation. Thus, the type of $\delta$ is a subsequence of $j_0,j_1,\ldots,j_{k-1}$, which is impossible by the first part of the claim and $(\clubsuit)$.
%
\end{proof}

\begin{clm} \label{clm:4}
If $k\geq 2$, then for $s, t\in \{1,\ldots,k-1\}$
\begin{itemize}
\item[(1)] if $s \ne t$ and $|s-t|\neq 1$, then $a^+(j_s,j_t) = 0$,
\item[(2)] if $s\neq 1$ and $t\neq k-1$, then $a^+(j_0,j_s) = a^+(j_k,j_t) = a^+(j_0,j_k) = 0$.
\end{itemize}
\end{clm} 
\begin{proof}[Proof of Claim \ref{clm:4}]
The claim follows by $(\clubsuit)$.
%
\end{proof}

The next claim follows immediately from Proposition \ref{prop:autf}.

\begin{clm} \label{clm:5}
Suppose $i\in I^+\setminus \{j_0,j_1,\ldots,j_k\}$ and $j\in I^+$, $j\ne i$ are such that $j$ is on the chain in the tree $G(a^+)$ joining $i$ and some color from $\{j_0,j_1,\ldots,j_k\}$ (i.e. the edge $\overrightarrow{(i,j)}$ is directed towards $\{j_0,j_1,\ldots,j_k\}$). Then $a^+(i,j)=1$.
\end{clm}

By Claim \ref{clm:5} there is a unique subtree $A'$ of $A$ corresponding to the code ${a^+}_{|\{j_0,\ldots,j_k\}^2}$. Hence, $A'$ is $\autf^+(A)$-invariant. What is left is to show that $A'$ is a biregular tree (which is clear when $n,m\geq 3$). If $m=n=2$, then $A'$ is a two-way infinite chain. Thus $\autf^+(A)$ leaves invariant a pair $\{b,b'\}$ of ends of $A$. As $\stab^{\autf^+(A)}(b)$ has index at most 2 in $\autf^+(A)$, it is a normal subgroup of $\autf^+(A)$. Hence $\autf^+(A)$ fixes $b$ (as $\autf^+(A)$ is a simple group), which is impossible. If $n=2$ and $m\geq 3$, then $A'$ is the $m$-regular tree.
\end{proof}

\begin{remark}
Notice that under the notation of Theorem \ref{thm:main}, $\stab^{\autf^+(A)}(A')$ is a normal subgroup of $\autf^+(A)$, so is trivial. Hence, Claim 6 from the proof can be strengthened: for $i\in I^+\setminus \{j_0,j_1,\ldots,j_k\}$ and $j\in I^+$, $j\ne i$, if $a^+(i,j)\ne 0$, then $a^+(i,j)=a^+(j,j)=1$.
\end{remark}

\section{Boundedly simple action on trees}

In this section we extend our results to some other groups acting on trees.

For a group $G$ acting on a tree $A$ we may consider the following coloring function \[f^{G} \colon \Som(A) \to \{\text{orbits of }G\text{ on }\Som(A)\}.\] The function $f^{G}$ is normal. 

\begin{definition} \label{def:k}
Denote by $\K$ the following class of groups: $G\in\K$ if and only if $G$ and every subgroup of $G$ of index 2 are virtually boundedly generated by finitely many conjugacy classes, that is $H=\prod_{1\leq i\leq n}h_i^H\cdot H_0$ for some finitely many $h_1,\ldots,h_n\in H$ and some finite $H_0\subset H$, where $H=G$ or $H\lhd G$ and $[G:H]=2$.
\end{definition}

The proof of the proposition below is standard.

\begin{proposition} \label{prop:k}
\begin{enumerate}
\item The class $\K$ contains all boundedly simple groups and all finite groups.
\item If $G\in\K$, then every quotient of $G$ is in $\K$.
\item If $G_1$ is a normal subgroup of $G$ and $G_1,G/G_1\in\K$, then $G\in\K$; that is $\K$ is closed under extensions.
\end{enumerate}
In particular $\K$ is closed under finite (semi)direct products.
\end{proposition}

\begin{theorem} \label{thm:act}
Suppose $A$ is a tree, $G<\aut(A)$, $G$ is in the class $\K$, $G$ leaves no nonempty proper subtree of $A$ invariant and does not stabilize any end of $A$. Then $A$ a uniform subdivision of an $(n,m)$-regular tree, for some $n,m\geq 3$.
\end{theorem}
\begin{proof} We may assume that $G$ is infinite.

\begin{claim}
$G\subseteq \aut_{f^G}^+(A)$ or $G\cap \aut_{f^G}^+(A)$ is a subgroup of $G$ of index 2.
\end{claim}
\begin{proof}[Proof of the claim]
Suppose that $G\not\subseteq \aut_{f^G}^+(A)$. Then \[\G_0 = G/\left(G\cap \aut_{f^G}^+(A)\right) = G/\aut_{f^G}^+(A)\] is in $\K$ (by Proposition \ref{prop:k}$(2)$) and $\G_0$ is a nontrivial subgroup of $\G = \aut_{f^G}(A)/\aut_{f^G}^+(A)$. By the main theorem of \cite{tits}, the group $\G$ is a free product $*_{i\in I}G_i$, for some index set $I$, where each $G_i$ is isomorphic to $\Z$ or to $\Z_2=\Z/2\Z$. By the Kurosh subgroup theorem, the group $\G_0$ is a free product of the form \[F(X)*\left(*_{j\in J}g_j^{-1}H_jg_j\right),\] where $F(X)$ is the free group, freely generated by $X\subseteq \G$, $J$ is some index set, $g_j\in \G$ and each $H_j$ is a subgroup of some $G_i$. It is enough to prove that either $\G_0\cong \Z_2$ or $\Z$ is a homomorphic image of $\G_0$ (notice that $\Z\not\in\K$). By the universal property of free products there is an epimorphism $A*B \to A\times B$, hence there exists an epimorphism \[\G_0 \to F(X)\times \prod_{j\in J}g_j^{-1}H_jg_j.\] If $X\ne\emptyset$ or some $H_j\cong\Z$, then the proof is finished. Otherwise, $\G_0\cong\Z_2$ or $\G_0$ has as a homomorphic image the group $\Z_2*\Z_2$. However, notice that $\Z_2*\Z_2\not\in\K$, as $\Z$ is a subgroup of index 2 in $\Z_2*\Z_2$.
\end{proof}

Let $H = G\cap \aut_{f^G}^+(A)$. As $H<\aut_{f^G}^+(A)$ and $G\in\K$, for some $N\in\N$ every translation from $H$ is the product of $N$ elements from \[S := \bigcup_{\{x,y\}\in\Are(A)} \stab^{\aut_{f^G}(A)}(x,y).\]

Since $\aut_{f^G}^+(A)$ is nontrivial (because $G$ is infinite), it contains some nontrivial rotation. By Lemma \ref{lem:rami}$(1)$, there is a ramification point $r\in\Som(A)$. Take an arbitrary $\alpha\in\aut_{f^G}^+(A)$. There is $h\in G$ with \[\alpha(r) = h(r).\] Hence (by Proposition \ref{prop:ram}), for some rotation $\beta\in\aut_{f^G}^+(A)$ fixing $r$, $\alpha = \beta\circ h$. The element $h = \beta^{-1}\alpha\in \aut_{f^G}^+(A)\cap G = H$ is a product of at most $N$ elements from $S$. Thus $\alpha$ is a composition of at most $N+2$ elements from $S$. Since $\aut_{f^{G}}^+(A)$ has property $(P)$, the assumption and Lemma \ref{lem:rami}$(1)$ imply that every translation from $\aut_{f^{G}}^+(A)$ is the product of $M$ elements each of which fixes pointwise a half-tree, for some $M\in\N$. The group $G$ leaves invariant no proper subtree and does not stabilize any end, so by \cite[Lemma 4.4]{tits} the same is true for $H$ and for $\aut_{f^{G}}^+(A)$. By Lemma \ref{lem:rami}$(2)$, $\aut_{f^{G}}^+(A)$ is boundedly simple. It is enough to apply now Theorem \ref{thm:main}.
\end{proof}

\begin{corollary} \label{cor:act}
Suppose $G$ is a group acting by automorphisms on a tree $A$ and $G\in\K$. Then $G$ fixes some vertex of $A$, or stabilizes some end of $A$, or the smallest nonempty $G$-invariant subtree $A' \subseteq A$ is a uniform subdivision of a biregular tree.
\end{corollary}
\begin{proof} By \cite[Corollary 3.5]{tits}, there is a nonempty minimal $G$-invariant subtree $A'$ of $A$, so Theorem \ref{thm:act} can be applied to $G/\stab^G(A')$ and $A'$.
\end{proof}



\begin{acknowledgements} 
I would like to thank Dugald Macpherson for useful conversations, Alexey Muranov for many suggestions for improving the presentation of the paper, especially for improving the bound in Theorem \ref{thm:bireg}, and Cong Chen for suggesting Lemma \ref{lem:cong}. 
\end{acknowledgements}

\end{document}